\documentclass[a4paper,12pt]{article} 
\newdimen\paperhight
\newcommand{\hf}{\frac{1}{2}}
\newcommand{\vv}{\langle -v,v\rangle}
\newcommand{\vu}{\langle -v,u\rangle}
\newcommand{\uu}{\langle -u,u\rangle}

\newcommand{\qed}{\begin{flushright}{\bf Q.E.D.}\end{flushright}}
\newcommand{\pr}{\par \vspace{3mm} \noindent {\bf [Proof]} \qquad}
\newcommand{\prend}{\hfill \qed}

\newcommand{\1}{{\bf 1}} 
\newcommand{\bw}{{\bf \omega}} 
\newcommand{\al}{\alpha} 
\newcommand{\be}{\beta} 
\newcommand{\de}{\delta}

\newcommand{\ga}{\gamma} 
\newcommand{\la}{\lambda}

\newcommand{\FP}{{\mathfrak p}} 
 
\newcommand{\CH}{{\cal H}}

\newcommand{\C}{{\mathbb C}} 
\newcommand{\Z}{{\mathbb Z}} 
\newcommand{\N}{{\mathbb N}}

\newcommand{\End}{{\rm End}}
 
\newcommand{\tr}{{\rm tr}}
 
\newcommand{\ch}{{\rm ch}} 

\newtheorem{thm}{Theorem}[section]
\newtheorem{prn}{Proposition}[section]
\newtheorem{dfn}{Definition}
\newtheorem{lmm}{Lemma}[section]

\usepackage{amssymb}
\usepackage[usenames]{color}
\usepackage{graphicx}

\begin{document}
\title{A modular invariance on the theta functions 
defined on vertex operator algebras} 
\author{Masahiko Miyamoto
\footnote{Supported by the Grants-in-Aids for Scientific Research, 
No. 09440004 and No. 10974001, The Ministry of Education, Science and 
Culture, Japan.}} 
\date{\begin{tabular}{c}
Institute of Mathematics \cr
University of Tsukuba \cr
Tsukuba 305, Japan \cr
\end{tabular} }
\maketitle
\begin{abstract}
In this paper, we introduce theta-functions of VOA-modules and 
show that the space spanned by them has a modular invariance property.   
If a VOA is 
a lattice VOA $V_L$ associated with an even lattice $L$ (cf. \cite{FLM}), 
then the above 
theta functions coincide to the ordinary theta functions 
$\theta(\tau, \vec{v})=\sum_{m\in L+\be}e^{\pi \imath <m,m>\tau+<m,\vec{v}>}$ 
for $<\be, L>\subseteq \Z$, $\vec{v}\in {\C}L$ and $\tau\in \CH$.  
\end{abstract}

\section{Introduction}
Throughout this paper, $V$ denotes a vertex operator 
algebra $(\oplus_{n=0}^{\infty} V_n,Y,\1,\bw)$ 
with central charge $c$ and $Y(v,z)=\sum v(n)z^{-n-1}$ denotes a vertex 
operator of $v$. (Abusing the notation, we will also use it 
for vertex operators of $v$ for $V$-modules.) 
$o(v)$ denotes the grade-keeping operator of $v$, 
which is given by $v(m-1)$ for $v\in V_m$ and defined by 
extending it for all elements of $V$ linearly.
In particular, $o(\bw)$ equals to $L(0)=\bw(1)$ for 
Virasoro element $\bw$ of $V$ and 
$o(v)=v(0)$ for $v\in V_1$. \hfill \\ 
In order to simplify the situation, we assume that $\dim V_0=1$ so that 
there is a constant $\langle v,u\rangle \in \C$ 
such that $v_1u=-\langle v,u\rangle\1$ for $v, u\in V_1$. 

We call $V$ a rational vertex operator algebra in case each $V$-module 
is a direct sum of simple modules. 
Define $C_2(V)$ to be the subspace of $V$ spanned by elements $u(-2)v$ 
for $u,v\in V$.  We say that $V$ satisfies {\it Condition $C_2$} if 
$C_2(V)$ has finite codimension in $V$.  
For a $V$-module $M$ with grading $M=\oplus M_m$, we define the 
formal character as 
\begin{equation}
\ch_qM=q^{-c/24}\sum \dim M_mq^m=\tr_M q^{-c/24+L(0)}.  
\end{equation}
In this paper, we consider these functions less formally by 
taking $q$ to be the usual local parameter $q=q_{\tau}=e^{2\pi \imath \tau}$ 
at infinity in the upper half-plane 
$$ \CH=\{\tau\in \C|\Im \tau>0\}.  $$

Although it is often said that a vertex operator algebra (VOA) 
is a conformal field theory with 
mathematically rigorous axioms, the axioms of VOA do not 
assume the modular invariance.  However, Zhu \cite{Z} showed the
modular ($SL_2(\Z)$) invariance of the space 
\begin{equation}\label{zhu}
<q_1^{|a_1|}\cdots q_n^{|a_n|}\tr_W Y(a_1,q_1)
\cdots Y(a_n,q_n)q^{L(0)-c/24} :\ W \mbox{ irr. $V$-mod}>
\end{equation}
for a rational VOA $V$ with central charge $c$ and $a_i\in V_{|a_i|}$
under Condition $C_2$, which 
are satisfied by many known examples, where $q_j=q_{z_j}=e^{2\pi \imath z_j}$ 
and $|a_i|$ denotes the weight of $a_j$. 
For example, the space 
$$<\ch_q W :\ W \mbox{  irreducible $V$-modules}> $$
is $SL_2(\Z)$-invariant. 
Recently, Dong, Li and Mason 
extended the Zhu's idea and proved a modular invariance of the space 
\begin{equation}
<\tr_U \phi^i q^{L(0)-c/24}: i\in \Z, \ U \mbox{ $\phi$-twisted modules}>
\end{equation}
by introducing the concept of $\phi$-twisted 
modules for a finite automorphism $\phi$, see \cite{DLiM}.  
An easiest example of 
automorphism of VOA is given by a vector $v\in V_1$ as 
$\phi=e^{2\pi \imath v(0)}$.  Especially, 
if the eigenvalues of $o(v)$ $(=v(0))$ on modules are in ${1\over n}\Z$, then 
the order of $e^{2\pi \imath o(v)}$ is finite.  So for a $V$-module $W$ and 
$u, v\in V_1$, we will define 
\begin{equation}
 Z_W(v;u;\tau)=\tr_W e^{2\pi\imath (o(v)-{\langle v,u\rangle \over 2})}
q^{L(0)+o(u)-(c+12\langle u,u\rangle )/24} 
\end{equation}
and we call $Z_W(v;0;\tau)\eta(\tau)^c$ a theta-function of $W$, where 
$u_1u=-\langle u,u\rangle\1$ and \\
$\eta(\tau)=q^{1/24}\prod_{n=1}^{\infty}(1-q^n)$ is 
the Dedekind $\eta$-function. 
It is worth noting that 
$c+12\langle u,u\rangle $ is the central 
charge of conformal element $\bw+L(-1)u$ and 
$o(\bw+L(-1)u)$ is equal to $L(0)+o(u)$, see \cite{DLnM}.

For example, 
let $V$ is a lattice VOA $V_{2\Z x}$ constructed 
from a $1$-dimensional lattice $L=2\Z x$ with $\langle x,x\rangle=1$. 
It has exactly four irreducible modules \cite{D}: 
$$W_0=V_{2\Z x}, W_1=V_{(2\Z+{1\over 2})x}, 
W_2=V_{(2\Z+1)x}, 
W_3=V_{(2\Z-{1\over 2})x}. $$
Let $\theta_{h,k}(z,\tau)$ $\left(:=\sum_{n\in \Z}\exp(\pi \imath(n+k)^2\tau
+2\pi(n+k)(z+k))\right)$ be theta functions for $h,k=0,\hf$. 
By the construction of lattice VOA (see \cite{FLM}), 
it is easy to check 
\begin{equation}
 \theta_{h,k}(z,\tau)
=\eta(\tau)((\imath)^{4hk}Z_{W_{2h}}(zx(-1)\1;0;\tau)
+(-1)^k(\imath)^{4hk}Z_{W_{2+2h}}(zx(-1)\1;0;\tau)) 
\end{equation}
for $h,k=0,\hf$ and  
their modular transformations 
\begin{equation}
\theta_{h,k}(z/\tau,-1/\tau)
=(i)^{4hk}(-\imath \tau)^{\hf}\exp(\pi \imath z^2/\tau)\theta_{k,h}(z,\tau) 
\end{equation}
are well known, (see \cite{M}).
In particular, there are constants $A_{k}^h\in \C$ such that 
\begin{equation}
 Z_{W_h}(zx(-1)\1;0;-1/\tau)
=\sum A_k^h Z_{W_h}(0;zx(-1)+\hf z^2,\tau). 
\end{equation}
Namely, the modular transformations of 
$Z_{W_h}(zx(-1)\1;0:\tau)$ are expressed by linear 
combinations of $Z_{W_k}(u;v;\tau)$ of (ordinary) modules $W_k$, but not 
twisted modules. 
By this result, for an automorphism $\phi=e^{v(0)}$, we can expect to 
obtain a modular transformation by using only the ordinary modules, 
which offers 
some information about twisted modules.  
This is the motivation of this paper and we will actually show that the above 
result is generally true, that is, 
we will prove that the following modular transformation by using 
Zhu's result (2). \hfill \\

\vspace{5mm}

\noindent
{\bf Main Theorem} \qquad 
Let $V$ be a rational vertex operator algebra with 
the irreducible modules $\{W_i:i=1,...,m \}$ and 
$u,v\in V_1$. 
Assume $v(0)v=v(0)u=u(0)v=u(0)u=0$ and $v(1)v,v(1)u,u(1)u\in {\C}\1$. 
If $V$ satisfies Zhu's finite condition $C_2$, 
then 
\begin{equation}
\{ Z_{W_h}(v;u;\tau): h=1,...,m \}
\end{equation}
satisfy a modular invariance, i.e., for 
$\al=\pmatrix{a&b\cr f&d}\in SL_2(\Z)$, there are constants 
$A_{\al,k}^h$ (see Theorem 4.1) such that 
\begin{equation}
 Z_{W_h}(v;u;{a\tau+b\over f\tau+d})
=\sum_{k=1}^m A^h_{\al,k}Z_{W_k}(av+bu;fv+du;\tau).
\end{equation}

\section{VOA}
A {\it vertex operator algebra} is a $\Z$-graded vector space:
\begin{equation}
  V=\oplus_{n\in \Z}V_n 
\end{equation}
satisfying $\dim V_n<\infty$ for all $n$ and $V_n=0 \mbox{  for  } n<<0$,  
equipped with a linear map 
\begin{eqnarray*}
&V \to &(\End V)[[z,z^{-1}]]  \hfill \\
&v \to &Y(v,z)=\sum_{n\in \Z}v(n)z^{-n-1} 
\end{eqnarray*} 
and with two distinguished vectors, {\it vacuum element} $\1\in V_0$ and 
{\it conformal vector} $\bw\in V_2$ satisfying 
the following conditions for $u,v\in V$: \hfill \\
\begin{eqnarray*}
&&u(n)v=0 \quad \mbox{ for $n$ sufficiently large}; \hfill \\
&&Y(\1,z)=1; \hfill \\
&&Y(v,z)\1\in V[[z]] \mbox{  and  }\lim_{z\to 0}Y(v,z)\1=v; \hfill \\
&&(z-x)^NY(v,z)Y(u,x)=(z-x)^NY(u,x)Y(v,z)  \quad 
\mbox{ for $N$ sufficiently large}, 
\end{eqnarray*} 
where $(z_1-z_2)^n\ (n\in \Z)$ are to be expanded in nonnegative integral 
powers of second variable $z_2$;
\begin{equation}
 [L(m),L(n)]=(m-n)L(m+n)+{1\over 12}(m3-m)\delta_{m+n,0}c  
\end{equation}
for $m,n\in \Z$, where $L(m)=\omega(m+1)$ and $c$ is called 
{\it central charge}; 
\begin{eqnarray*} 
&&L(0)v=nv \mbox{  for  }v\in V_n;  \hfill \\
&&{d\over dz}Y(v,z)=Y(L(-1)v,z).  
\end{eqnarray*}
This completes the definition. \hfill \\

We also have the notion of modules: \hfill \\
Let $(V,Y,\1,\omega)$ be a vertex operator algebra. A {\it weak} module $W$ of 
$(V,Y,\1,\omega)$ is a $\C$-graded vector space: 
\begin{equation}
 W=\oplus_{n\in \C}W_n   
\end{equation}
equipped with a linear map 
\begin{eqnarray*}
&V \to &(\End(W))[[z,z^{-1}]] \hfill \\
&v \to &Y^W(v,z)=\sum_{n\in \Z}v^W(n)z^{-n-1} \qquad (v_n\in End(W)) 
\end{eqnarray*} 
satisfying the following conditions for $u,v\in V$ and $w\in W$: \hfill \\
For $v\in V, w\in M$, $v^W(m)w=0$ for $m>>0$. \hfill \\
$$ Y^W(\1,z)=1; $$
$$ L^W(0)w=nw \mbox{  for  }w\in W_n, L^W(0)=\omega^W(1);  $$ 
$$ {d\over dz}Y^W(v,z)=Y(L(-1)v,z)  $$   
and the following Jacobi idenity holds.
\begin{eqnarray*}
&&z_0^{-1}\delta({z_1-z_2\over z_0})Y^W(u,z_1)Y^W(v,z_2)
-z_0^{-1}\delta({z_2-z_1\over -z_0})Y^W(v,z_2)Y^W(u,z_1) \hfill \\
&&=z_2^{-1}\delta({z_1-z_0\over z_2})Y^W(Y(u,z_0)v,z_2). 
\end{eqnarray*}

A weak module $W$ is called a {\it module} if 
every finitely generated weak submodule $M=\oplus_{r\in \C} M_r$ of $W$ 
satisfies \hfill \\
(1)  $\dim M_r<\infty$,  \hfill \\
(2)  $M_{r+n}=0 \mbox{  for $n\in \Z$ sufficiently large, }$ \hfill \\
for any $r\in \C$.

\section{Formal power series}
We use the notation $q$ and $q_z$ to denote $e^{2\pi \imath \tau}$ and 
$e^{2\pi \imath z}$, 
respectively.

In this paper, the formal power series 
$$ P_2(q_z,q)=
(2\pi \imath )^2\sum_{n=1}^{\infty}({nq_z^n \over 1-q^n}
+{nq_z^{-n}q^n\over 1-q^n})$$
plays an essential role,  
where ${1\over 1-q^n}$ is understood as $\sum_{i=0}^{\infty}q^i$. 
The limit of $P_2(q_z,q)$ (which we still denote as $P_2(q_z,q)$) relates to 
$\FP(z,\tau)$ by 
$$ P_2(q_z,q)=\FP(z,\tau)+G_2(\tau),   $$
where 
$$G_2(\tau)={\pi^2\over 2}+\sum_{m\in \Z-\{0\}}\sum_{n\in \Z}
{1\over (m\tau+n)^2}$$ 
is the Eisenstein series and 
$\FP(z,\tau)$ is the Weierstrass $\FP$-function 
$$ \FP(z,\tau)={1\over z^2}+\sum_{(m,n)\not=(0,0)}
({1\over (z-m\tau-n)^2}-{1\over (m\tau+n)^2}).   $$
It is known that 
$$ G_2({a\tau+b\over f\tau+d})=(f\tau+d)^2G_2(\tau)-2\pi \imath f(f\tau+d)$$ 
and 
$$ \FP({z\over f\tau+d},{a\tau+b\over f\tau+d})=(f\tau+d)^2\FP(z,\tau).$$
In particular, 
\begin{equation} P_2({z\over f\tau+d},{a\tau+b\over f\tau+d})
=(f\tau+d)^2P_2(z,\tau)-2\pi \imath f(f\tau+d).  
\end{equation}

In this paper, we will use variables $\{z_1,...,z_n\}$ 
and calculate the products of formal power series $P_2(q_{z_i-z_j},\tau)$. 
In order to simplify notation, we will use 
a transposition $(i,j)$ of symmetric groups $\Sigma_n$ on 
$\Omega=\{1,...,n\}$.  For $\{(i_{11},i_{12}),...,(i_{t1},i_{t2})\}$ with 
$i_{s1}<i_{s2}$ and $i_{ab}\not=i_{cd}$ for $(a,b)\not=(c,d)$, 
we view $\sigma=(i_{11},i_{12})\cdots (i_{t1},i_{t2})$ as an 
involution (element of order $2$) of $\Sigma_n$ and 
denote $\prod_{j=1}^t P_2(q_{z_{t2}-z_{t1}},\tau)$ by 
$\prod_{i<\sigma(i)} P_2(q_{z_{\sigma(i)}-z_i},\tau)$.  
Let $I(n)$ denote the set of all elements $g$ in $\Sigma_n$ with $g^2=1$.  
For $\sigma\in \Sigma_n$, set
 \begin{eqnarray}
&m(\sigma)=\{i\in \Omega| \sigma(i)\not=i\} \hfill \\
&f(\sigma)=\{i\in \Omega|\sigma(i)=i\}. 
\end{eqnarray} 
For $\sigma_1, \sigma_2\in \Sigma_n$, $\sigma_1, \sigma_2$ are 
called to be {\it disjoint} if 
$m(\sigma_1)\cap m(\sigma_2)=\emptyset$. 
$\sigma_1+\cdots+\sigma_n=\sigma$ expresses that 
$\{\sigma_1,...,\sigma_n\}$ are mutually disjoint and 
the product $\sigma_1\cdots\sigma_n$ is equal to $\sigma$.

For $v\in V_1$, $v(0)$ acts on the finite dimensional homogeneous 
subspaces $V_m$ and satisfies $[v(0), v(m)]=(v(0)v)(m)$ and 
\begin{eqnarray*}
&v(0)\bw&=v(0)\bw(-1)\1=-[\bw(-1),v(0)]\1
=-\sum_{i=0}^{\infty}(-1)^i(\bw(i)v)(-1-\imath )\1 \hfill \\
&&=-(\bw(0)v)(-1)\1+(\bw(1)v)(-2)\1=-v(-2)\1+(w(1)v)(-2)\1=0. 
\end{eqnarray*}   
Therefore, 
$e^{v(0)}=\sum_{n=0}^{\infty} {v(0)^n \over n!}$ is well defined 
and satifies that  
$e^{v(0)}\bw=\bw$ and 
$e^{v(0)}(u_mw)=(e^{v(0)}u)_m(e^{v(0)}w)$. In particular, 
$e^{v(0)}$ is an automorphism of $V$.  

\begin{dfn}
For a $V$-module $W$ and $u,v\in V_1$, define 
\begin{equation}
Z_W(v;u;\tau)
=\tr_W e^{2\pi \imath (v(0)-\langle v,u\rangle/2)}
q^{(u(0)-\langle u,u\rangle/2)+L(0)-c/24}. 
\end{equation}
We set $\theta_W(v,\tau)=Z_W(v;0;\tau)\eta(\tau)^c$ and call it 
a theta-function of $W$, where \\
$\eta(\tau)=q^{1/24}\prod_{n=1}^{\infty}(1-q^n)$ is the 
Dedekind $\eta$-function. 
\end{dfn}

For example, 
let $V$ be a lattice VOA $V_{2\Z x}$ associated with 
a $1$-dimensional lattice $L=2\Z x$ with $\langle x,x\rangle=1$, then 
$$W_0=V_{2\Z x}, W_1=V_{(2\Z+{1\over 2})x}, 
W_2=V_{(2\Z+1)x}, 
W_3=V_{(\Z-{1\over 2})x}$$
are the irreducible $V_{2\Z x}$-modules by \cite{D} and 
$\langle x(-1)\1, x(-1)\1\rangle=-1$. 
We then have 
\begin{eqnarray*}
&\theta_{0,0}(z,\tau)&
=\sum_{n\in \Z}\exp(\pi \imath n^2\tau+2\pi \imath nz)\hfill \\
&&=\theta_{W_0}(zx(-1)\1,\tau)+\theta_{W_2}(zx(-1)\1,\tau)) \hfill \\
&\theta_{0,\hf}(z,\tau)&
=\sum\exp(\pi \imath n^2\tau+2\pi \imath n(z+\hf))\hfill \\
&&=\theta_{W_0}(zx(-1)\1,\tau)-\theta_{W_2}(zx(-1)\1,\tau)) \hfill \\
&\theta_{\hf,0}(z,\tau)&=\sum\exp(\pi \imath (n+\hf)^2\tau+2\pi \imath 
(n+\hf)z)\hfill \\
&&=\theta_{W_1}(zx(-1)\1,\tau)+\theta_{W_3}(zx(-1)\1,\tau)) \hfill \\
&\theta_{\hf,\hf}(z,\tau)&=\sum\exp(\pi \imath (n+\hf)^2\tau+2\pi \imath 
(n+\hf)(z+\hf)) \hfill \\
&&=\imath\theta_{W_1}(zx(-1)\1,\tau)-\imath \theta_{W_3}(zx(-1)\1,\tau)) 
\end{eqnarray*} 
and their modular transformations 
\begin{eqnarray*}
&&\theta_{0,0}(z/\tau,-1/\tau)
=(-\imath \tau)^{\hf}\exp(\pi \imath z^2/\tau)\theta_{0,0}(z,\tau) \hfill \\
&&\theta_{0,\hf}(z/\tau,-1/\tau)
=(-\imath \tau)^{\hf}\exp(\pi \imath z^2/\tau)\theta_{1,0}(z,\tau) \hfill \\
&&\theta_{\hf,0}(z/\tau,-1/\tau)
=(-\imath \tau)^{\hf}\exp(\pi \imath z^2/\tau)\theta_{0,1}(z,\tau) \hfill \\
&&\theta_{\hf,\hf}(z/\tau,-1/\tau)
=-(-\imath \tau)^{\hf}\exp(\pi \imath z^2/\tau)\theta_{1,1}(z,\tau)  
\end{eqnarray*} 
are well known, (see \cite{M}). 
Therefore, the modular transformations 
$Z_{W_h}(zx(-1)\1;0:{-1\over \tau})$ are expressed by linear 
combinations of $Z_{W_k}(u;v;\tau)$ of (ordinary) modules $W_k$, but not 
twisted modules.

\section{Modular invariance}
In this section, we will prove a modular invariance. 
Throughout this section, we assume : \hfill \\
(A1)  $V=\oplus_{n=0}^{\infty}V_n$ is a rational VOA; \hfill \\
(A2)  $\{W^1,...,W^m\}$ is the set of all irreducible $V$-modules;  \hfill \\
(A3)  fix $v_1,...,v_n\in V_1$ satisfying $v_r(0)v_j=0$ and 
$v_r(1)v_j\in {\C}\1$ for any $r,j$. \hfill \\

By (A3), we have 
$$[v_r(m),v_j(n)]=\sum {m\choose k}(v_r(k)v_j)(m+n-\imath )
=m(v_r(1)v_j)(m+n-1)
=\delta_{m,-n}m\langle -v_r,v_j\rangle, $$
where $v_r(1)v_j=\langle -v_r,v_i\rangle\1$.   

For a grade-keeping endomorphism formal power series 
$\psi\in \End(W)((q_{z_1},...,q_{z_n},q_{\tau}))$ of $W$, set
$$S_W(\psi;z_1,...,z_n,\tau)=\tr_W\psi 
Y(v_1,q_{z_1})\cdots Y(v_n,q_{z_n})q_{z_1}\cdots 
q_{z_n}q^{L(0)-{c\over 24}}.   $$

By the same argument as in the proof of Proposition 4.3.2 in \cite{Z}, 
we have the following: 

\begin{prn}
Assume $[\psi, v_r(n)]=0$ for $r=1,...,n$. Then we have : 
\begin{equation}
S_W(\psi;z_1,z_2,...,z_n,\tau)
=\!\sum_{\sigma\in I(n)}\!S_W\!\left(\psi\! \prod_{r\in f(\sigma)}
\!o(v_r)\prod_{j<\sigma(j)}\!
\left(\langle -v_j,v_{\sigma(j)}\rangle 
{P_2(z_{\sigma(j)}z_j^{-1},\tau)\over (2\pi\imath)^2}\right);\tau\right),
\end{equation}
where $I(n)$ is the set of all elements $\sigma$ in the symmetric group 
$\Sigma_n$ on $n$ point set $\Omega=\{1,2,...,n\}$ with $\sigma^2=1$ 
and $f(\sigma)$ denotes the 
set of fixed points of $\sigma$.   
\end{prn}

\pr   For $k\in \Z$, we have 
\begin{eqnarray*}
& &S_W(\psi v_1(k)q_{z_1}^{-k};z_2,...,z_n,\tau) \hfill \\
&=&\!\!\tr_W  \psi v_1(k)q_{z_1}^{-k}Y(v_2,q_{z_2})\cdots 
Y(v_n,q_{z_n})q_{z_1}\cdots q_{z_n}q^{L(0)-c/24} \hfill \\
&=&\!\!\tr_W  \psi[v_1(k),Y(v_2,q_{z_2})\cdots Y(v_n,q_{z_n})]q_{z_1}^{-k}q_{z_2}
\cdots q_{z_n}q^{L(0)-c/24} \hfill \\
& &+\tr_W  \psi Y(v_2,q_{z_2})\cdots Y(v_n,q_{z_n})v(k)q_{z_1}^{-k}q_{z_2}
\cdots q_{z_n}q^{L(0)-c/24} \hfill \\
&=&\!\!\sum_{j=2}^n\sum_{i\in \N}{k\choose i}q_{z_j}^{k-\imath }
\tr_W  \psi Y(v_2,q_{z_2})\cdots Y(v_1(i)v_j,q_{z_j})\cdots 
Y(v_n,q_{z_n})q_{z_1}^{-k}q_{z_2}\cdots q_{z_n}q^{L(0)-c/24} \hfill \\
&&+\tr_W \psi Y(v_2,q_{z_2})\cdots Y(v_n,q_{z_n})q_{z_1}^{-k}q_{z_2}
\cdots q_{z_n}q^{L(0)-c/24}v(k)q^{k} \hfill \\
&=&\!\!\sum_{j=2}^n kq_{z_j}^{k-1}q_{z_1}^{-k}\tr_W \psi Y(v_2,q_{z_2})\cdots 
Y(\langle -v_j,v_1\rangle \1,q_{z_j})\cdots Y(v_n,q_{z_n})q_{z_2}\cdots 
q_{z_n}q^{L(0)-c/24} \hfill \\
&&+\tr_W  \psi  v_1(k)q_{z_1}^{-k}Y(v_2,q_{z_2})\cdots Y(v_n,q_{z_n})q_{z_2}
\cdots q_{z_n}q^{L(0)-c/24}q^k \hfill \\
&=&\!\!\sum_{j=2}^n\tr_W \langle -v_j,v_1\rangle kq_{z_j-z_1}^k \psi 
Y(v_2,q_{z_2})\cdots \widehat{Y(v_j,q_{z_{j}})}
\cdots Y(v_n,q_{z_n})q_{z_2}...\widehat{q_{z_j}}...
q_{z_n}q^{L(0)-c/24} \hfill \\
&&+\tr_W  \psi  v_1(k)q_{z_1}^{-k}Y(v_2,q_{z_2})\cdots Y(v_n,q_{z_n})q_{z_2}
\cdots q_{z_n}q^{L(0)-c/24}q^k \hfill \\
&=&\!\!\sum_{j=2}^n\tr_W \langle -v_j,v_1\rangle kq_{z_j-z_1}^k 
\psi Y(v_2,q_{z_2})\cdots 
\widehat{Y(v_j,q_{z_{j}})}\cdots Y(v_n,q_{z_n})
q_{z_2}\cdots \widehat{q_{z_{j}}}\cdots q_{z_n}q^{L(0)-c/24} \hfill \\
&&+\tr_W  \psi  v_1(k)q_{z_1}^{-k}Y(v,q_{z_2})\cdots Y(v,q_{z_n})q_{z_2}
\cdots q_{z_n}q^{L(0)-c/24}q^k \hfill \\
&=&\!\!\sum_{j=2}^n \langle -v_1,v_j\rangle kq_{z_j-z_1}^k
S(\psi;z_2,...,\widehat{z_{j}},...,z_n,\tau)
+S(\psi v_1(k)q_{z_1}^{-k};z_2,...,z_n,\tau)q^k,  
\end{eqnarray*} 
where $\widehat{q_{z_j}}$ means that we take off the term $q_{z_j}$. 

Hence, for $k\not=0$, we have:   
\begin{equation}
S_W(\psi v_1(k)q_{z_1}^{-k};z_2,...,z_n,\tau)
=\sum_{j=2}^n\langle -v_j,v_1\rangle {kq_{z_j-z_1}^{k} \over 1-q^k}
S_W(\psi;z_2,...,z_{j-1},z_{j+1},...,z_n,\tau). 
\end{equation}

Therefore, we have:
\begin{eqnarray*}
&&\quad S_W(\psi;z_1,z_2,...,z_n,\tau)\hfill \\
&=&S_W(\psi v_1(0);z_2,...,z_n,\tau)
+\sum_{k\not=0}S_W(\psi v(k)q_{z_1}^{-k};z_2,...,z_n,\tau)\hfill \\
&=&S_W(\psi v_1(0);z_2,...,z_n,\tau)
+\sum_{k\not=0}\sum_{j=2}^n
\langle -v_j,v_1\rangle {kq_{z_j-z_1}^{k} \over 1-q^k}
S_W(\psi;z_2,...,z_{j-1},z_{j+1},...,z_n,\tau)  \hfill \\
&=&S_W(\psi v_1(0);z_2,...,z_n,\tau)
\!-\!\sum_{j=2}^n\langle v_u,v_1\rangle 
{P_2(q_{z_j-z_1},\tau) \over (2\pi\imath)^2} 
S_W(\psi;z_2,...,z_{j-1},z_{j+1},...,z_n,\tau)  
\end{eqnarray*} 

By substituting $\psi v_1(0)$ into $\phi$ 
and repeating these steps, we have:
\begin{equation}
S_W(\psi;z_1,z_2,...,z_n,\tau)
=\sum_{\sigma\in I(n)}
\prod_{i<\sigma(i)}\left( \langle -v_i,v_{\sigma(i)}\rangle 
{P_2(q_{z_{\sigma(i)}-z_i},\tau) \over 
(2\pi\imath)^2}\right)
S_W(\psi \prod_{i\in f(\sigma)}v_i(0),\tau). 
\end{equation}
\prend

The main result we quote from \cite{Z} is 

\begin{thm}[Zhu]
Let $(V,Y,\1,\bw)$ be a rational VOA satisfying (A1) and (A2). Then 
for any $\al=\pmatrix{a &b \cr f& d}\in SL_2(\Z)$, we have  
\begin{equation}
S_{W_h}(1,{z_1\over f\tau+d},...,{z_n\over f\tau+d},{a\tau+b\over 
f\tau+d})=(f\tau+d)^n\sum_{j=1}^m A^h_{\al,k}S_{W_k}(1;z_1,...,z_n,\tau)  
\end{equation} 
where the $A^h_{\al,k}$ are constants depending only on $\al,h,k$. 
\end{thm}

To simplify notation, $S_k(\cdots)$, $\al(v)$ and 
$d(\al)$ denote $S_{W_k}(\cdots)$, $(f\tau+d)o(v)$ and $f\tau+d$, 
respectively.  
Set  
\begin{eqnarray*}
&D(r,j)=&\langle -v_r,v_j\rangle({1\over 2\pi\imath})^2
(d(\al))^2P_2(q_{z_r-z_j},\tau), \hfill \\
&E(r,j)=&\langle -v_r,v_j\rangle ((d(\al))^2P_2(q_{z_r-z_j},\tau)
-({1\over 2\pi \imath})f(d(\al))), \hfill \\ 
&D_{\sigma}=&\prod_{j<\sigma(j)}D(\sigma(j),j), \mbox{  and } \hfill \\
&E_{\sigma}=&\prod_{j<\sigma(j)}E(\sigma(j),j).
\end{eqnarray*}

\begin{lmm}
If $|m(\sigma)|=2p$, then 
\begin{equation}
\sum_{\sigma_1+...+\sigma_t=\sigma}(-1)^{t}
E_{\sigma_t}\cdots E_{\sigma_2}E_{\sigma_1}
=(-1)^pE_{\sigma}.
\end{equation}
\end{lmm}

\pr We first note that 
$E_{\sigma_t}\cdots E_{\sigma_2}E_{\sigma_1}=E_{\sigma}$. 
Therefore, we have to count the number of $E_{\sigma}$ in the left side.  
We will prove it by induction on $p$.  If $p=1$, then it is trivial. 
For $\sigma=(r_1,r_2)\cdots(r_{2p-1},r_{2p})$, the number of $\sigma_1$ 
with $\sigma_1+\cdots+\sigma_t=\sigma$ and $|m(\sigma_1)|=2r$ 
is ${p\choose r}$. 
Therefore, by induction, we have:
\begin{equation}
\sum_{\sigma_1+...+\sigma_t=\sigma}(-1)^{t}
E_{\sigma_t}\cdots E_{\sigma_2}E_{\sigma_1}
=-\sum_{j=1}^p{p\choose j}(-1)^{p-j}E_{\sigma}
=-(-(-1)^p)E_{\sigma}=(-1)^pE_{\sigma}. 
\end{equation}
\prend

For $\al=\pmatrix{a&b\cr f&d}\in SL_2(\Z)$, we have 

\begin{lmm}  
\begin{equation}
S_h(\prod_{r=1}^n o(v_r); {a\tau+b\over f\tau+d})
=\sum_{k=1}^m A_{\al,k}^h \sum_{\sigma\in I(n)} \prod_{j<\sigma(j)}
({f(d(\al))\over 2\pi\imath })\langle -v_j,v_{\sigma(j)}\rangle 
S_k(\prod_{r\in f(\sigma)}\al(v_r);\tau).
\end{equation}
\end{lmm}

\pr
By Theorem 3.1 and Proposition 3.1, we have:
\begin{eqnarray*}
&&\sum_{\sigma\in I(n)}\left(\prod_{r<\sigma(r)}
\left(\langle -v_{\sigma(r)},v_r\rangle({1\over 2\pi\imath })^2 
P_2(q_{{z_{\sigma(r)}-z_r \over 
f\tau+d}}, {a\tau+b\over f\tau+d})\right)\right)
S_h(\prod_{j\in f(\sigma)}o(v_j),{a\tau+b\over f\tau+d}) \hfill \\
&=&S_h(1;{z_1\over f\tau+d},...,{z_n\over f\tau+d},{a\tau+b\over f\tau+d})\hfill \\
&=&(d(\al))^n\sum_{k=1}^mA^h_{\al,k}S_k(1;z_1,...,z_n,\tau) \hfill \\
&=&(d(\al))^n\sum_{k=1}^mA^h_{\al,k}\sum_{\sigma\in I(n)}
\left(\prod_{r<\sigma(r)}\left(\langle -v_r,v_{\sigma(r)}\rangle 
({1\over 2\pi\imath })^2P_2(q_{z_{\sigma(r)}-z_r},\tau)\right)
\right)S_k(\prod_{j\in f(\sigma)}o(v_j);\tau). 
\end{eqnarray*} 

Hence, we have
\begin{eqnarray*}
&&S_h(\prod_{r=1}^n o(v_r); {a\tau+b\over f\tau+d}) \hfill \\
&=&\sum_{k=1}^mA^h_{\al,k}\sum_{\sigma\in I(n)}\left(\prod_{r<\sigma(r)}
\left(\langle -v_r,v_{\sigma(r)}\rangle (d(\al))^2
({1\over 2\pi\imath })^2 P_2(q_{z_{\sigma(r)}-z_r},\tau)\right)\right)
S_k(\prod_{j\in f(\sigma)}\al(v_j)),\tau) \hfill \\
&&-\sum_{1\not=\sigma\in I(n)}\left(\prod_{r<\sigma(r)}
\left(\langle -v_r,v_{\sigma(r)}\rangle 
((d(\al))^2({1\over 2\pi\imath })^2 P_2(q_{z_{\sigma(r)}-z_r},\tau)-
({1\over 2\pi\imath }) f(d(\al)))\right)\right)
\times  \hfill \\
 &&\times S_h(\prod_{j\in f(\sigma)} o(v_j);{a\tau+b\over f\tau+d}) \hfill \\
&=&\sum_{k=1}^m A^h_{\al,k}\sum_{\sigma\in I(n)}\left(
\prod_{r<\sigma(r)}D(\sigma(r),r)\right)S_k(\prod_{j\in f(\sigma)}\al(v_j), 
\tau)  \hfill \\
&&-\sum_{1\not=\sigma\in I(n)}\left(\prod_{r<\sigma(r)}E(\sigma(r),r)\right)
S_h(\prod_{j\in f(\sigma)} o(v_j); {a\tau+b\over f\tau+d}) \hfill \\ 
&=&\sum_{k=1}^m A^h_{\al,k}\sum_{\sigma\in I(n)}D_{\sigma}
S_k(\prod_{j\in f(\sigma)}\al(v_j); \tau)
-\sum_{1\not=\sigma\in I(n)}E_{\sigma}
S_h(\prod_{r\in f(\sigma)} o(v_r); {a\tau+b\over f\tau+d}).
\end{eqnarray*}

Substituting the above equation into the last term, we have 
\begin{eqnarray*}
 &&S_h(\prod_{r=1}^n o(v_r); {a\tau+b\over f\tau+d}) \hfill \\
&=&\sum_{k=1}^m A^h_{\al,k}\sum_{\sigma\in I(n)}D_{\sigma}
S_k(\prod_{j\in f(\sigma)}\al(v_j);\tau)  \hfill \\
&&-\sum_{1\not=\sigma\in I(n)}E_{\sigma}
\sum_{k=1}^m A^h_{\al,k}
\sum_{\sigma'\in I(n), m(\sigma')\cap m(\sigma)=\emptyset}
D_{\sigma'}S_k(\prod_{r\in f(\sigma'+\sigma)}\al(v_r); \tau)  \hfill \\
&&+\sum_{1\not=\sigma\in I(n)}E_{\sigma}
\sum_{1\not=\sigma'\in I(n),m(\sigma')\cap m(\sigma)=\emptyset}
E_{\sigma'}S_h(\prod_{j\in f(\sigma)\cap f(\sigma')}o(v_j);
{a\tau+b\over f\tau+d}) \}. 
\end{eqnarray*} 

Repeating the above steps and by the equations 
\begin{eqnarray*}
&&S_h(o(v);{a\tau+b\over f\tau+d})=\sum A^h_{\al,k}S_k(o(v);\tau)\hfill \\
&&S_h(1;{a\tau+b\over f\tau+d})=\sum A^h_{\al,k}S_k(1;\tau).
\end{eqnarray*} 
in \cite{Z}, we have 
\begin{eqnarray*}
&&S_h(\prod_{r=1}^n o(v_r);{a\tau+b\over f\tau+d})  \hfill \\
&&=\sum_{k=1}^m A^h_{\al,k}
S_k(\sum_{\sigma\in I(n)}D_{\sigma}\prod_{j\in f(\sigma)}\al(v_j);\tau)  
\hfill \\
&&\quad -\sum_{1\not=\sigma\in I(n)}E_{\sigma}
 S_h(\prod_{r\in f(\sigma)} o(v_r);{a\tau+b\over
f\tau+d})  \hfill \\
&&=\sum_{k=1}^m A^h_{\al,k}
\{\sum_{\sigma_1,...,\sigma_t\in I(n):|f(\sigma_2+...\sigma_t)|\leq n-2}(-1)^{t-1}
E_{\sigma_t}\cdots E_{\sigma_2}D_{\sigma_1}
S_k(\prod_{r\in \cap_{j=1}^t f(\sigma_j)}\al(v_r);\tau)  \hfill \\
&&\mbox{(if $n$ is odd)  } -\sum_{j=1}^n
\sum_{\sigma_1,...,\sigma_t\in I(n):f(\sigma_1+...+\sigma_t)
=\{j\}}(-1)^{t}
E_{\sigma_t}\cdots E_{\sigma_2}E_{\sigma_1}S_k(\al(v_j);\tau)\}  \hfill \\
&&\mbox{(if $n$ is even)  }
-\sum_{\sigma_1,...,\sigma_t\in I(n):f(\sigma_1+...+\sigma_t)=\emptyset}
(-1)^{t}
E_{\sigma_t}\cdots E_{\sigma_2}E_{\sigma_1}S_k(1;\tau)\} \hfill \\
&&\mbox{by Lemma 3.1}  \hfill \\
&&=\sum_{k=1}^m A_{\al,k}^h\{ \sum_{\sigma\in I(n)}
\sum_{\sigma_1,\sigma_2\in I(n):\sigma_1+\sigma_2=\sigma, |f(\sigma_1)|\geq 2}
(-1)^{|\sigma_1|/2}E_{\sigma_1}D_{\sigma_2}
S_k(\prod_{r\in f(\sigma)}\al(v_r);\tau)  \hfill \\
&&\mbox{(if $n$ is odd)  }
-\sum_{j=1}^n\sum_{\sigma\in I(n): f(\sigma)=\{j\}}(-1)^{(n-1)/2}E_{\sigma}
S_k(\al(v_j);\tau)\}   \hfill \\
&&\mbox{(if $n$ is even)  }
-\sum_{\sigma\in I(n): f(\sigma)=\emptyset}(-1)^{n/2}E_{\sigma}S_k(1;\tau)\}  
\hfill \\
&&=\sum_{k=1}^m A_{\al,k}^h\{\sum_{\sigma\in I(n)} 
\prod_{j<\sigma(j)}({1\over 2\pi \imath})
(f(d(\al))\langle -v_j,v_{\sigma(j)}\rangle )^p
S_k(\prod_{r\in f(\sigma)}\al(v_r);\tau)\}
\end{eqnarray*} 
\prend

\noindent
{\bf Theorem A}\qquad  For $\al=\pmatrix{a&b\cr f&d}\in SL_2(\Z)$, we have 
\begin{equation}
 S_h(e^{o(v)};{a\tau+b\over f\tau+d})
=\sum_{k=1}^m A^h_{\al,k}S_k(e^{{\langle -v,v\rangle\over 2}
({1\over 2\pi \imath})f(d(\al))+(d(\al))o(v)}, \tau).
\end{equation}
In particular, 
\begin{equation}
 Z_{W_h}(v;0;{a\tau+b\over f\tau+d})=\sum_{k=1}^mA_{\al,k}^h 
Z_{W_k}(dv;fv;\tau). 
\end{equation}

\pr 
Let's calculate the coefficient of 
$(\langle -u,v \rangle({1\over 2\pi \imath}) f(d(\al)))^p(\al(v))^r$ 
in (4.2) by setting $n=2p+r$ and $v_1=\cdots=v_{n}=v$.

The number of involutions $\sigma$ with $|f(\sigma)|=r$ is 
${2p+r\choose r}\times {(2p+r)!\over p!2^p}$.  
Therefore, we have 
\begin{eqnarray*}
&&S_h(e^{o(v)};{a\tau+b\over f\tau+d})
=S_h(\sum_{n=0}^{\infty}{1\over n!}\prod_{i=1}^n o(v);{a\tau+b\over f\tau+d})
=\sum_{n=0}^{\infty}{1\over n!}S_h(\prod_{i=1}^n o(v);
{a\tau+b\over f\tau+d})  \hfill \\
&&=\sum_{r,p\in \N}{1\over n!}\sum_{k=1}^mA^h_{\al,k}
\sum_{\sigma\in I(n)}(\langle -v,v\rangle 
{f(d(\al))\over 2\pi\imath})^{|m(\sigma)|/2}S_k(
(\al(v))^{|f(\sigma)|};\tau)  \hfill \\
&&=\sum_{k=1}^mA^h_{\al,k}
\sum_{r,p\in \N}{1\over (r+2p)!}{r+2p\choose r}{(2p)!\over p!2^p}
(\langle -v,v \rangle {f(d(\al))\over 2\pi\imath})^{|m(\sigma)|/2}
S_k( (\al(v))^{|f(\sigma)|};\tau)  \hfill \\
&&=\sum_{k=1}^mA^h_{\al,k}
\sum_{r,p\in \N}{1\over (r+p)!}{r+p\choose p}
({\langle -v,v \rangle \over 2}{f(d(\al))\over 2\pi\imath})^{p}
S_k((\al(v))^{r};\tau)  \hfill \\
&&=\sum_{k=1}^mA^h_{\al,k}
\sum_{n=0}^{\infty}\sum_{r=0}^n
{1\over n!}{n\choose p}
({\langle -v,v\rangle\over 2} 
{f(d(\al))\over 2\pi\imath})^{n-r}S_k((\al(v))^{r};\tau)  \hfill \\
&&=\sum_{r=1}^mA^h_{\al,k}S_k(\sum_{n=0}^{\infty}
{1\over n!}(({\langle -v,v\rangle \over 2}{f(d(\al))\over 2\pi\imath})
+\al(v))^n;\tau)  \hfill \\
&&=\sum_{r=1}^mA^h_{\al,k}S_k(e^{{\langle -v,v\rangle \over 2}
{f(d(\al))\over 2\pi\imath}+\al(v)};\tau), 
\end{eqnarray*} 
since 
\begin{equation}
{1\over (r+2p)!}{r+2p\choose r}{(2p)!\over p!2^p}={1\over (r+2p)!}
{(r+2p)!\over k!(2p)!}{(2p)!\over p!2^p}={1\over (r+p)!}{r+p\choose r}
{1\over 2^p}.
\end{equation}
\prend

Let's show the final version of our result. \hfill \\ 

\noindent
{\bf Main Theorem}\qquad 
For $\al=\pmatrix{a&b\cr f&d}\in SL_2(\Z)$, we have 
\begin{equation}
Z_{W_h}(v;u;{a\tau+b\over f\tau+d})=\sum_{k=1}^m A_{\al,k}^h
Z_{W_k}(dv+bu;fv+au;\tau),  
\end{equation}
where $A_{\al,k}^h$ are the coefficients in the equations given by 
Zhu \cite{Z}. \hfill \\

\pr
Fix $s,t\in \N$ so that $n=s+t$.   
Assume that 
$$v_1=v_2=...=v_{s},\ (\mbox{ say, }=v), \quad 
u_{s+1}=\cdots =v_{s+t},\ (\mbox{ say, }=u). $$
For $p$, $q$, $r\in \Z$, set $k=s-2p-r\geq 0$, $h=t-2q-r\geq 0$. 

In order to simplify notation, for $u,v\in V_1$ and $\sigma\in I(n)$, 
we will use the following notation: 
\begin{eqnarray*} 
&&\alpha(-v,u)={f(d(\al))\over 2\pi \imath}
\langle -v,u\rangle, \hfill \\
&&\al(v)=(d(\al))o(v) \hfill \\
&&m11\sigma=|\{r\in \Omega| r<\sigma(r)\leq s\}|, \hfill \\
&&m12\sigma=|\{r\in \Omega| r\leq s<\sigma(r)\}|, \hfill \\
&&m22\sigma=|\{r\in \Omega| s<r<\sigma(i)\}|, \hfill \\
&&f1\sigma=|\{r\in \Omega| r=\sigma(r)\leq s\}| \mbox{  and} \hfill \\
&&f2\sigma=|\{r\in \Omega| s<i=\sigma(r)\}|.  
\end{eqnarray*} 

Then by lemma 4.2, we have
\begin{eqnarray*}
&&\mbox{}\quad S_h(o(v)^{s}o(u)^{t};{a\tau+b\over f\tau+d})\hfill \\
&&=\sum_{g=1}^m A_{\al,g}^h \sum_{\sigma\in I(s+t)}
\alpha(-v,v)^{m11\sigma}
\alpha(-v,u)^{m12\sigma}
\alpha(-u,u)^{m22\sigma}
S_g(\al(v)^{f1\sigma}\al(u)^{f2\sigma};\tau).  
\end{eqnarray*} 
Hence, 
\begin{eqnarray*}
&&\mbox{}\quad S_h(o(v)^{s}o(u)^{t};{a\tau+b\over f\tau+d})  \hfill \\
&&=\sum_{g=1}^m A_{\al,g}^h 
\sum_{r,p,q}{(2p+r+k)!(2q+r+h)! \over  k!h!r!p!q!2^{p+q}}
\alpha(-v,v)^{p}
\alpha(-v,u)^{r}
\alpha(-u,u)^{q}
S_g(\al(v)^{k}\al(u)^{h};\tau).   
\end{eqnarray*}

Expanding the exponential, we have:
\begin{eqnarray*}
&&S_h(e^{2\pi \imath(o(v)+\langle -v,u\rangle/2)
+2\pi\imath\tau(o(u)+\langle -u,u\rangle/2)};\tau) \hfill \\
&&=\tr_{W_h}\sum_{n=0}^{\infty}{(2\pi\imath)^n\over n!}
(o(v)+\tau o(u)+{\vu \over 2}+\tau{\uu \over 2})^n
q^{L(0)-c/24}  \hfill \\
&&=\tr_{W_h}\sum_{\al,\be,\ga,\delta}
{(2\pi\imath)^{\al+\be+\ga+\de}\over \al!\be!\ga!\de!}
o(v)^{\al}(\tau o(u))^{\be}({\vu \over 2})^{\ga}
({\uu \over 2}\tau)^{\delta}q^{L(0)-c/24}. 
\end{eqnarray*} 
Hence, for $\al={-1\over \tau}$, we have:
 \begin{eqnarray*}
&&S_h(e^{2\pi i(o(v)+\langle -v,u\rangle/2+\tau(o(u)+\uu /2))};{-1\over \tau}) \hfill \\
&&=\sum_{\al,\be,\ga,\delta}
{(2\pi \imath)^{\al+\be+\ga+\de}\over \al! \be! \ga! \delta!}
S_h(o(v)^{\al}({-1\over \tau} o(u))^{\be}
({\vu \over 2})^{\ga}({\uu\over 2})^{\delta}({-1\over \tau})^{\delta},
{-1\over \tau})  \hfill \\
&&=\sum_{\ga,\delta,\al,\be}
{(2\pi\imath)^{\al+\be+\ga+\de}\over \al!\be!\ga!\de!}
({\uu\over 2})^{\delta}( {-1\over \tau})^{\delta}({\vu\over 2} )^{\ga}
({-1\over \tau})^{\be}S_h(o(v)^{\al}o(u)^{\be};{-1\over \tau}) \hfill \\
&&=\sum_{\ga,\delta,\al,\be}
{(2\pi \imath)^{\al+\be+\ga+\de}\over \al!\be!\ga!\de!}
({\uu\over 2})^{\delta}( {-1\over \tau})^{\delta}({\vu\over 2} )^{\ga}
({-1\over \tau})^{\be}\sum_{g=1}^m A_{\al,g}^h
\sum_{p,q,r,h,k:2p+r+k=\al, 2q+r+h=\be} \hfill \\
&&{\al!\be!\over h!k!r!p!q!2^{p+q}}\times 
S_g((\tau{\vv\over 4\pi\imath} )^p(\tau{\vu\over 4\pi\imath} )^r
(\tau{\uu\over 4\pi\imath} )^q(\tau o(v))^k(\tau o(u))^h,\tau)  \hfill \\
&&=\sum_{\ga,\delta,p,q,r,h,k}{(2\pi \imath)^{(p+k+r+q+h+\ga+\de)}\over 
\ga!\de!h!k!r!p!q!}
({-1\over \tau}({\uu \over 2}))^{\delta}({\vu \over 2})^{\ga}\times  \hfill \\
&&\times \sum_{g=1}^m A_{\al,g}^hS_g(
(\tau{\vv\over 2} )^p(-\vu)^r({1\over \tau}{\uu\over 2})^q
(\tau o(v))^k(-o(u))^h;\tau)  \hfill \\
&&=\sum_{g=1}^m A_{\al,g}^h S_g(
e^{2\pi\imath({-1\over \tau}\uu/2 +\vu /2+
\tau\vv/2 -\vu +{1\over \tau}\uu/2 +\tau o(v)-o(u)) };\tau)  \hfill \\
&&=\sum_{g=1}^m A_{\al,g}^h S_g(
e^{2\pi\imath(-\vu /2+
\tau\vv/2 +\tau o(v)-o(u))};\tau)  \hfill \\
&&=\sum_{g=1}^m A_{\al,g}^h S_g(
e^{2\pi\imath(o(-u)+\langle u,v\rangle/2)+2\pi\imath\tau(o(v)+\vv/2)};\tau)  
\end{eqnarray*} 
Namely, 
\begin{equation}
Z_{W_h}(v;u;{-1\over \tau})=
\sum A_{\al,k}^h Z_{W_k}(-u;v;\tau).  
\end{equation}
On the other hand, 
\begin{equation}
 Z_{W_h}(v;u;\tau+1)=Z_{W_h}(v+u;u;\tau).  
\end{equation}
Therefore, for $\al=\pmatrix{a&b\cr f&d}$, we have 
\begin{equation} 
Z_{W_h}(v;u;{a\tau+b\over f\tau+d})=\sum_{k=1}^mA_{\al,k}^h Z_{W_k}(
dv+bu;fv+au;\tau). 
\end{equation}
\prend

\noindent
{\bf Example} 
\qquad Let $V$ be a lattice VOA $V_{2\Z x}$ with $\langle x,x\rangle=1$ 
and $W$ a module $V_{\Z x}$.  
For $z\in \C$, set
$$zx=zx(-1)\1\in (V_{2\Z x})_1. $$
Then it is easy to check by the definition of lattice VOAs that 
$${\tr} o(\1)e^{2\pi izx(0)}q^{L(0)-1/24}|_M
={1\over \eta(\tau)}\theta(\tau,z). $$
For $\al=\pmatrix{0&-1 \cr 1&0}$, 
$${1\over \eta({-1\over \tau})}\theta({-1\over \tau},z)=\sum_M 
\la_M {\tr} 
o({\1})q^{L(0)+zx(0)-{1\over 24}+{z^2\over 2}}|_M, $$
where $M$ runs over the following modules 
$$\{V_L, V_{x+L}, V_{\hf x+L}, V_{-\hf x+L} \}.$$
Taking $\tau=i$ and several $z$, we have the known formula:
$${1\over \eta({-1\over \tau})}\theta({-1\over \tau},z)
={1\over \eta(\tau)}\theta(\tau,\tau z)e^{\pi iz^2\tau}. $$


\begin{thebibliography}{99}
\bibitem[B]{B}
R.~E.~Borcherds, 
{\it Vertex algebras, Kac-Moody algebras, and the Monster}, \hfill \\
Proc.~Natl.~Acad.~Sci. USA 83 (1986), 3068-3071.

\bibitem[CN]{CN}
J.~H.~Conway and S.~P.~Norton, {\it Monstrous moonshine,} 
Bull.~London Math. Soc. 11(1979), 308-339.

\bibitem[D]{D}
C.~Dong, 
{\it Vertex algebras associated with even lattices}, 
J. Algebra 160 (1993), 245-265. 

\bibitem[DLiM]{DLiM}
C.~Dong, H.~Li and G.~Mason, 
{\it Modular-invariance of trace functions in orbifold theory}, 
preprint.


\bibitem[DLnM]{DLnM}
C.~Dong, Z.~Lin and G.~Mason,
"On vertex operator algebras as $sl_2$-modules, in: 
{\it Groups, Difference Sets, and the Monster, Proc. of a Special Research 
Quarter at The Ohio State University, Spring 1993}, ed. by K.T.~Arasu, 
J.F.~Dillon, K.~Harada, S.~Sehgal and R.~Solomon, Walter de Gruyter, 
Berlin-New York, 1996, 349-362. 



\bibitem[FHL]{FHL}
I.~B.~Frenkel, Y.-Z.~Huang and J.~Lepowsky, "On axiomatic approaches to
vertex operator algebras and modules",  
{\it Memoirs Amer. Math. Soc.} {\bf 104}, 1993.

\bibitem[FLM]{FLM}
I.~B.~Frenkel, J.~Lepowsky, and A.~Meurman,
"Vertex Operator Algebras and the Monster", 
Pure and Applied Math., Vol. 134, 
Academic Press, 1988.

\bibitem[M]{M}
D.~Mumford, 
"Tata Lectures on Theta I", 
Progress in Mathematics, Vol. 28, Birkh\"{a}user, 1983. 

\bibitem[Z]{Z}
Y.~Zhu, {\it Modular invariance of characters of vertex operator algebras}, 
J. Amer. Math. Soc. 9 (1996), 237-302. 
\end{thebibliography}
\end{document}